\documentclass[12pt]{article}
\usepackage{amssymb}

\newcommand{\eps}{\varepsilon}
\renewcommand{\phi}{\varphi}

\newtheorem{theorem}{Theorem}

\newtheorem{lemma}{Lemma}

\title{On rich lines in grids}

\author{Evan Borenstein\thanks{Summer funding supported by an NSF VIGRE
grant.} and Ernie Croot \thanks{Supported in part by 
an NSF grant.}}

\begin{document}

\maketitle

\section{Introduction}

In \cite{erdos_szemeredi}, Erd\H os and Szemer\'edi proved the following
result, which has led to a remarkable number of profound developments
in the field of additive combinatorics:

\begin{theorem}  There is some absolute constant $\eps > 0$ such that 
if $A$ is a set of real numbers, $|A| \geq 2$, then either the sumset
$A+A$ or the product set $A.A$, has size at least $|A|^{1+\eps}$.
\end{theorem}

In \cite{elekes}, Elekes gave a brilliantly elegant proof of this theorem
using the Szemer\'edi-Trotter incidence theorem \cite{szemeredi}, and
was able to show that 
$$
|A+A| \cdot |A.A|\ \gg\ |A|^{5/2},
$$
from which it follows that
$$
\max(|A+A|,\ |A.A|)\ \gg\ |A|^{5/4}.
$$
The key fact that Elekes needed for his proof, and which is a weak corollary
of the Szemer\'edi-Trotter incidence theorem, at least as far as 
just getting a non-trivial bound of the sort
$$
|A+A| \cdot |A.A|\ \gg\ |A|^{2+\eps},
$$
is the following basic claim.
\bigskip

\noindent {\bf Claim 1.}  There are absolute constants $\eps > 0$ and
$\delta > 0$ such that if $A$ and $B$ are sets of $n$ real numbers,
and $n$ is sufficiently large (in terms of $\eps$ and $\delta$), 
then any set of at least $n^{2-\eps}$ distinct 
lines contains a member that hits the grid in fewer than $n^{1-\delta}$
points.  In other words, one cannot have a collection of $n^{2-\eps}$ lines
whereby all are ``$n^{1-\delta}$-rich'' in the grid $A \times B$.
\bigskip

Actually, Elekes's proof only needs the following even weaker claim.
\bigskip

\noindent {\bf Claim 2.}  There exist absolute constants $\eps > 0$
and $\delta > 0$ so that the following holds for all integers $n$
sufficiently large:  Suppose that $A$ and $B$ are sets of real numbers
of size $n$, and that one has a family of lines such that
\bigskip

$\bullet$ There are at least $n^{1-\eps}$ distinct slopes among them; and,

$\bullet$ every line is parallel to at least $n^{1-\eps}$ others.
\bigskip

\noindent Then, at least one of the lines must hit the grid $A \times B$
in fewer than $n^{1-\delta}$ points.  In other words, not all the
lines can be $n^{1-\delta}$-rich in the grid.
\bigskip  

In the present paper we prove the following theorem, which shows that
it is possible to considerably strengthen this second claim;
futhermore, our theorem is not the sort that is quickly deducible
from the Szemer\'edi-Trotter incidence theorem:

\begin{theorem} \label{main_theorem}  
For every $\eps > 0$, there exists $\delta > 0$
so that the following holds for all $n$ sufficiently large:  Suppose
that $A$ and $B$ are sets of real numbers of size $n$, and that one 
has a family of lines such that
\bigskip

$\bullet$  There are at least $n^\eps$ distinct slopes among them; and,

$\bullet$  every line is parallel to at least $n^\eps$ others.
\bigskip

\noindent Then, at least one of the lines must hit the grid
$A \times B$ in fewer than $n^{1-\delta}$ points.  
\end{theorem}

Our theorem is related to a conjecture of Solymosi 
(see \cite[Conj. 3.10]{elekes2} for details), which we modify and 
extend to make it better fit the context of the above results.
\bigskip

\noindent {\bf Solymosi's Conjecture.}  For every $\eps > 0$,
there exists $\delta > 0$, such that the following holds for all integers
$n$ sufficiently large:  Suppose $A$ and $B$ are sets of real numbers 
of size $n$, and suppose that one has a collection of $n^\eps$ lines 
in general position (that is, no pair is parallel, and no three meet 
at a point).  Then, not all of the lines can be $n^{1-\delta}$-rich 
in the grid $A \times B$.
\bigskip

This conjecture of Solymosi easily implies our main theorem 
(Theorem \ref{main_theorem}) above, for if one has a family of lines as
described by our theorem, then it is a simple matter to select one 
line from each of $\gg n^{\eps/3}$ groups of parallel lines in such a
way that one produces a collection in general position
(first, select a single line of slope $\lambda_1$; then, select a line
of slope $\lambda_2 \neq \lambda_1$; then, select a line of slope 
$\lambda_3 \not \in \{\lambda_1,\lambda_2\}$ such that the three lines
do not have a common intersection point; then, select a line of slope
$\lambda_4 \not \in \{\lambda_1,\lambda_2,\lambda_3\}$...). 

\subsection{Remarks}

Our proof makes use of several standard methods in additive combinatorics,
though is quite intricate and technical.  In particular, some of our 
approaches are similar to those appearing in the well-known paper of
Bourgain, Katz and Tao \cite{BKT}, as was pointed out to us by 
P. M. Wood.  Even so, we do not assume any results more sophisticated 
than the Szemer\'edi-Trotter theorem.  It was pointed out to us recently 
by T. Tao that perhaps we could make use of a particular sum-product 
ideas of Bourgain to give a simpler proof; however, we decided to 
present here our original approach.

It is possible that perhaps some of the ideas of Harald Helfgott 
\cite{helfgott} might allow us to give a shorter proof, as part of our
argument can be phrased in terms of growth and generation in subgroups
of $GL_2({\mathbb R})$.     
 
\section{Proof of the main theorem} 

The first step in our proof is to reduce from the case of working with
grids $A \times B$ to grids $A \times A$.  This is easily handled by
simply letting $C = A \cup B$, and then noting that the hypotheses of
our theorem imply that we have a family of rich lines passing through
the grid $C \times C$.  Upon rescaling $n$ to $|A \cup B| \leq 2n$, 
we see that we could have just assumed that our grid was $A \times A$
(or $C \times C$) all along. 

\subsection{Producing new rich lines from old ones} \label{combine_subsection}

In our proof we will be combining together lots of pairs of rich
lines, possibly of different slope:  Given a line $\ell$
hitting $A \times A$ in some points, we let 
\begin{eqnarray*}
X(\ell)\ &=&\ {\rm projection\ of\ } \ell \cap (A \times A)\ {\rm onto\
the\ x-axis}; \\
Y(\ell)\ &=&\ {\rm projection\ of\ } \ell \cap (A \times A)\ {\rm onto\
the\ y-axis}.
\end{eqnarray*}

If two lines
$$
\ell\ :\ y\ =\ \lambda x + \mu\ \ {\rm and\ \ }
\ell'\ :\ y\ =\ \lambda' x + \mu',
$$
have the property that
$$
|Y(\ell) \cap Y(\ell')|\ =\ ``{\rm large}",
$$
then there will be lots of triples
$$
(x,z,y)\ \in\ A \times A \times A
$$
satisfying
$$
\lambda x + \mu\ =\ y\ =\ \lambda' z + \mu'.
$$
So, the new line
$$
z\ =\ (\lambda/\lambda') x + (\mu - \mu')/\lambda'
$$
also hits the grid $A \times A$ in many points.
\bigskip

A convenient way of keeping track of the new rich lines that we
can produce from old ones is to use matrix notation:  We form the
association
$$
y\ =\ \lambda x + \mu\ \ \leftrightarrow\ \ \left [ \begin{array}{cc}
\lambda & \mu \\ 0 & 1 \end{array} \right ].
$$
Then, when we combine together lines as above, the new line we get
will be the one associated to a certain product of
matrices; specifically,
\begin{eqnarray*}
&& y\ =\ (\lambda/\lambda')x + (\mu - \mu')/\lambda'\\
&& \hskip1in \leftrightarrow\ \
\left [ \begin{array}{cc} \lambda/\lambda' & (\mu - \mu')/\lambda' \\
0 & 1 \end{array}\right ]\ =\
\left [ \begin{array}{cc} \lambda' & \mu' \\ 0 & 1 \end{array} \right ]^{-1}
\left [ \begin{array}{cc} \lambda & \mu \\ 0 & 1 \end{array} \right ]
\end{eqnarray*}
\bigskip

A basic fact, which is an easy consequence of the Cauchy-Schwarz 
inequality, is the following lemma:

\begin{lemma} \label{overlap_lemma}  Given lines
$$
\ell_1,\ ...,\ \ell_K,
$$
each hitting a grid
$$
A\ \times\ A,
$$
in at least
$$
n^{1 - \delta_0}\ {\rm points},
$$
we have that at least
$$
K^2 n^{-2\delta_0}/2
$$
of the pairs $(\ell_i,\ell_j)$ have the property that
\begin{equation} \label{pair_bound}
|Y(\ell_i) \cap Y(\ell_j)|\ \geq\ n^{1-2\delta_0}/2.
\end{equation}
\end{lemma}

If the lines $\ell_1,...,\ell_K$ have slopes $\lambda_1,...,\lambda_K$,
respectively, then upon combining it with our preceeding observations,
we deduce that there are lots of lines of slope 
$\lambda_i/\lambda_j$, for lots of pairs $(i,j)$,
such that each is at least $n^{1-2\delta_0}/2$ rich in the grid 
$A \times A$.  

\subsection{Passing to a set of rich lines with usable properties} 
\label{usable_subsection}

Given $\eps > 0$, we let $\delta' > 0$ denote some parameter that we
will choose later.  Then, given $\eps, \delta' >0$ we let $\delta > 0$
be some parameter chosen later.  We will show that if $\delta > 0$
is small enough, and if (as stated in the hypotheses of our theorem)
we had a set of lines involving $n^\eps$ slopes, each parallel to at 
least $n^\eps$ others, each $n^{1-\delta}$-rich in the grid, 
then in fact there would have to exist at least $n^4$ lines, 
each hitting $A \times A$ in at least two points.  This clearly cannot 
happen, because there are fewer lines hitting the grid in two points 
than there are ordered pairs of points of the grid; there are 
$n^2$ points of the grid, and therefore $n^4$ ordered pairs.  
This will prove our theorem.

So, we assume that $\eps > 0$ is given, and then we will select 
$\delta' > 0$ as small as needed, and then choose $\delta > 0$ even
smaller later.
\bigskip

We begin by letting $L_1(\lambda)$ denote the set of our lines having
slope $\lambda$.  We note that
$$
|L_1(\lambda)|\ \geq\ n^\eps,
$$
where $\lambda$ is one of the slopes of our set of lines.
To make certain later estimates easier, we will trim our list of
lines so that for each slope $\lambda$ we have 
$$
|L_1(\lambda)|\ =\ \lceil n^\eps \rceil.
$$
Denote our initial set of slopes by $\Lambda_1$.
\bigskip

Using Lemma \ref{overlap_lemma}, we can easily deduce that there are at 
least 
$$
|\Lambda_1|^2 n^{-O(\delta)}
$$
ordered pairs 
$$
(\lambda,\lambda')\ \in\ \Lambda_1 \times \Lambda_1,
$$
for which there are at least
$$
|L_1(\lambda)| \cdot |L_1(\lambda')| n^{-O(\delta)}\ \sim\ n^{2\eps - O(\delta)}
$$
pairs of lines 
\begin{equation} \label{ellpairs}
(\ell, \ell')\ \in\ L_1(\lambda) \times L_1(\lambda')
\end{equation}
satisfying
\begin{equation} \label{lambda'}
|Y(\ell) \cap Y(\ell')|\ \geq\ n^{1-O(\delta)}.
\end{equation}
Note that each of these intersections gives rise to a line having slope 
$\lambda/\lambda'$ that hits $A \times A$ in $n^{1-O(\delta)}$ points. 

When such a pair $(\lambda,\lambda')$ has the above property we will
say that it is ``good for step 1''.  Note that our definition of 
``good'' is dependent upon the implied constants in the big-ohs --
for our purposes, the implied constants in the ``good for step $i$'' can 
all be taken to be $1000^i$.

If a pair $(\lambda,\lambda')$ is good for step $1$, and if in addition
we have that the number of distinct lines of slope $\lambda/\lambda'$ 
produced by combining pairs $(\ell,\ell')$ satisfying 
(\ref{ellpairs}) and (\ref{lambda'}) is at least 
\begin{equation} \label{LLL}
n^{\eps (1 + \delta')},
\end{equation}
we will say that $(\lambda,\lambda')$ is ``very good for step 1''.   

Let us suppose that all but at least a fraction $n^{-O(\delta)}$
of the ``good'' pairs $(\lambda,\lambda')$ are, in fact, ``very good''.
Let $\Lambda_2$ denote these ``very good'' pairs, and note that we are
saying
$$
|\Lambda_2|\ \geq\ |\{ {\rm good\ pairs} \}| n^{-O(\delta)}\ \geq\ 
|\Lambda_1|^2 n^{-O(\delta)}.
$$
\bigskip

For $\theta \in \Lambda_2$, say $\theta = (\lambda,\lambda')$, we 
let $L_2(\theta)$ denote those lines produced by combining together 
pairs of lines, one from
$L_1(\lambda)$ and the other from $L_1(\lambda')$.  Note that for all
$\theta \in \Lambda_2$ we have, by (\ref{LLL}), that 
$$
|L_2(\theta)|\ \geq\ n^{\eps(1+\delta')}.
$$ 
And, as with the set of lines $L_1(\lambda)$, we trim our set of lines
(in an arbitrary manner) so that for every such $\theta$ we have that
$$
|L_2(\theta)|\ =\ \lceil n^{\eps(1+\delta')} \rceil.
$$

It is easily deduced from Lemma \ref{overlap_lemma} that there are at least
$$
|\Lambda_2|^2 n^{-O(\delta)}
$$
ordered pairs 
$$
(\theta,\theta')\ \in\ \Lambda_2 \times \Lambda_2,
$$
for which there are at least
$$
|L_2(\theta)| \cdot |L_2(\theta')| n^{-O(\delta)}\ \sim\ 
n^{2\eps(1 + \delta') - O(\delta)}
$$
pairs of lines 
$$
(\ell, \ell')\ \in\ L_2(\lambda) \times L_2(\lambda')
$$
satisfying
$$
|Y(\ell) \cap Y(\ell')|\ \geq\ n^{1-O(\delta)}.
$$

When such a pair $(\theta,\theta')$ has the above property we will say
that it is ``good for step 2'', and we say that it is 
``very good for step 2'' if the set of rich lines that it produces 
has size at least
$$
n^{\eps (1 + \delta')^2}.
$$

We will repeat the above process we have started as above, 
by defining $\Lambda_3$ to be the set of all ``very good for step 2'' pairs
$\beta = (\theta,\theta') \in \Lambda_2 \times \Lambda_2$, and we will let 
$L_3(\beta)$ be those lines produced by combining together ones from
$L_2(\theta)$ with $L_(\theta')$, and then trimming the list so that 
$$
|L_3(\beta)|\ =\ \lceil n^{\eps (1 + \delta')^2} \rceil.
$$ 
\bigskip

It is clear that we can continue the above process, producing sets 
$$
\Lambda_4,\ \Lambda_5,\ ...,\ {\rm where\ } \Lambda_i\ \subseteq\ \Lambda_{i-1}
\times \Lambda_{i-1},
$$
and sets 
$$
L_3(\alpha_3),\ L_4(\alpha_4),\ ...,\ {\rm where\ } \alpha_i\ \in\ \Lambda_i.
$$
However, the process cannot go on for too long, since we always have the
upper bound
$$
|L_t(\alpha)|\ \leq\ n^4,
$$
since the lines of $L_t(\alpha)$ will hit the grid in at least two points.
In fact, 
$$
t\ \ll\ T\ := (1/\delta') \log(4/\eps).
$$

Well, the above sequence of $\Lambda_j$'s and $L_j(\alpha_j)$'s is not
quite what we want, because for later arguments we will need that the
sequence terminates with $t > k$, for some $k = k(\eps)$ depending only
on $\eps > 0$.  The way we get around this is as follows:  Going back to
how our sequences of $\Lambda_j$'s and $L_j(\alpha_j)$'s are defined, if
we are willing to allow the $\Lambda_j$, $j=1,2,...,k$ to merely contain 
``good for step $j$'' pairs, instead of ``very good for step $j$'' pairs,
then the problem of stopping at time $t \leq k$ is avoided.  There is 
the issue of how to trim the sets $L_2(\alpha_2), ..., L_k(\alpha_k)$ in 
the right way.  To solve this problem, we merely trim them  
so that they each contain $n^{\eps - O(\delta)}$ lines, which is
easily guaranteed.  Furthermore, by choosing $\delta' > 0$ small enough, 
we can still have that for $j > k$ and $\theta \in \Lambda_j$,
$$
|L_j(\theta)|\ =\ \lceil n^{\eps (1 + \delta')^j} \rceil,
$$
the reason being that for small $\delta' > 0$, the $(1 + \delta')^k$
can be made as close to $1$ as needed. 
\bigskip

Before unraveling what this all means, we make one more observation:  
An element $\theta \in \Lambda_i$ corresponds to a pair of elements
of $\Lambda_{i-1}$, and each member of the pair itself corresponds 
to pairs of elements of $\Lambda_{i-2}$,
and so on; so, in the end, an element of $\theta \in \Lambda_i$ in fact
corresponds to a sequence of elements of $\Lambda_1$ of length 
$2^{i-1}$.  Say the sequence is 
$$
\lambda_1,\ ...,\ \lambda_{2^{i-1}}.
$$
Then, the lines it corresponds to all have slope
$$
\lambda_1 \cdots \lambda_{2^{i-2}}\ /\ \lambda_{2^{i-2}+1} \cdots 
\lambda_{2^{i-1}}.
$$ 
\bigskip

When our above process terminates at time $t$ satisfying 
$$
k\ <\ t\ \ll\ T,
$$
we will have that the following all hold:
\bigskip

$\bullet$  First, for at least 
$$
|\Lambda_1|^{2^{t-2}} n^{-O_t(\delta)}
$$
sequences
$$
\lambda_1,\ ...,\ \lambda_{2^{t-2}}\ \in\ \Lambda_1
$$
we will have a set of lines of slope 
$$
\lambda_1 \cdots \lambda_{2^{t-3}}\ /\ \lambda_{2^{t-3}+1} \cdots 
\lambda_{2^{t-2}}
$$
that are $n^{1-O_t(\delta)}$-rich in our grid $A \times A$.  
\bigskip

$\bullet$ Second, there are at least 
$$
|\Lambda_1|^{2^{t-1}} n^{-O_t(\delta)}
$$
pairs of sequences
$$
\lambda_1,\ ...,\ \lambda_{2^{t-1}} \in \Lambda_1,\ {\rm and\ } 
\lambda'_1,\ ...,\ \lambda'_{2^{t-1}} \in \Lambda_1,
$$
corresponding to a pair of elements 
$$
(\nu_1,\nu_2)\ \in\ \Lambda_{t-1} \times \Lambda_{t-1},
$$
that are ``good for step $t$'' but not ``very good for step $t$'' 
(since otherwise we could continue the iteration for another step).  
For such a pair, suppose that 
our $n^{1-O_t(\delta)}$-rich lines corresponding to $\nu_1$ are of the form
\begin{equation} \label{form1}
y\ =\ (\lambda_1 \cdots \lambda_{2^{t-3}} / \lambda_{2^{t-3}+1} \cdots 
\lambda_{2^{t-2}}) x\ +\ B_{\nu_1},
\end{equation}
and those corresponding to $\nu_2$ are of the form
\begin{equation} \label{form2}
y\ =\ (\lambda'_1 \cdots \lambda'_{2^{t-3}} / \lambda'_{2^{t-3}+1} 
\cdots \lambda'_{2^{t-2}}) x + B_{\nu_2}.
\end{equation}
Then, since the pair $(\nu_1,\nu_2)$ is ``good for step $t$'', we 
have that there are 
$$
|B_{\nu_1}| \cdot |B_{\nu_2}| n^{-O_t(\delta)}
$$
ordered pairs of lines, one corresponding to $\nu_1$ and the other
to $\nu_2$, such that when combined, give us an $n^{1-O_t(\delta)}$-rich
line of the form
$$
y\ =\ \alpha x + (b_1 - b_2)/\beta,
$$
where
$$
\alpha\ =\ \lambda_1 \cdots \lambda_{2^{t-3}} \lambda'_{2^{t-3}+1} \cdots 
\lambda'_{2^{t-2}} / \lambda'_1 \cdots \lambda'_{2^{t-3}} 
\lambda_{2^{t-3}+1} \cdots \lambda_{2^{t-2}},
$$
where 
$$
b_1\ \in\ B_{\nu_1},\ b_2\ \in\ B_{\nu_2},\ {\rm and\ where\ } 
\beta\ =\ \lambda'_1 \cdots \lambda'_{2^{t-3}} / \lambda'_{2^{t-3}+1} 
\cdots \lambda'_{2^{t-2}}.
$$
Furthermore, since the pair $(\nu_1,\nu_2)$ is not ``very good for
step $t$'', we have that the possibilities for the difference 
$b_1 - b_2$ is at most 
$$
n^{\eps (1 + \delta')^t}\ \leq\ |L_{t-1}(\nu_1)|^{1 + \delta'}\ =\ 
|B_{\nu_1}|^{1+\delta'}.
$$

What this means is that the ``additive energy'' between the sets
$B_{\nu_1}$ and $B_{\nu_2}$ must be ``large''.  In fact, because
there are so many pairs $(\nu_1,\nu_2)$, there must exist 
$\nu_1 \in \Lambda_{t-1}$ such that there are at least 
$$
|\Lambda_{t-1}| n^{-O_t(\delta)}
$$  
choices for $\nu_2 \in \Lambda_t$, such that we have the following
lower bound for the additive energy:
\begin{eqnarray*}
E(B_{\nu_1}, B_{\nu_2})\ &=&\ |\{ (b_1,b_2,b_3,b_4) \in 
B_{\nu_1} \times B_{\nu_1} \times B_{\nu_2} \times B_{\nu_2}
\ :\ b_1 - b_3\ =\ b_2 - b_4 \}| \\
&\geq&\ |B_{\nu_1}|^{3 - O(\delta')}. 
\end{eqnarray*}

We now require the following standard lemma.

\begin{lemma}  Suppose that $X$ and $Y$ are sets of size $M$, such that
$$
E(X,Y)\ =\ |\{ (x,x',y,y') \in X \times X \times Y \times Y\ :\ 
x - y = x' - y'\}|\ \geq\ c M^3.
$$
Then, there is some translate $u$ such that
$$
|(X+u) \cap Y|\ \geq\ c M.
$$
\end{lemma}

\noindent {\bf Proof of the Lemma.}  Another way of writing the additive
energy is 
$$
E(X,Y)\ =\ \sum_{u \in X \atop v \in Y} | (X - u) \cap (Y - v)|.
$$
So, by simple averaging, among the $M^2$ pairs 
$(u,v) \in X \times Y$, there exists one for which
$$
|(X - u + v) \cap Y|\ =\ |(X - u) \cap (Y - v)|\ \geq\ c M;
$$
\hfill $\blacksquare$
\bigskip

So, for some fixed $\nu_1 \in \Lambda_{t-1}$, and for 
$|\Lambda_{t-1}| n^{-O_t(\delta)}$ elements $\nu_2 \in \Lambda_{t-1}$, there
exist translates $\tau(\nu_2)$ for which
$$
|B_{\nu_1} \cap (B_{\nu_2} + \tau(\nu_2))|\ \geq\ 
|B_{\nu_1}| n^{-O_t(\delta')}.
$$ 
\bigskip

We now arrive at the following basic claim.
\bigskip

\noindent {\bf Claim 3.}  Under the hypotheses of our theorem, 
there are distinct slopes 
$$
\theta_1,...,\theta_N, 
$$
where 
$$
N\ >\ n^{\eps - O(\delta)},
$$
such that for 
$$
m\ =\ 2^{t-2},
$$
at least $N^{m-O(\delta)}$ of the $m$-fold products 
$\theta_{i_1}\cdots \theta_{i_m}$, we have a set of 
$n^{1-O(\delta)}$-rich lines of the form
$$
y\ =\ \theta_{i_1}\cdots \theta_{i_m} x\ +\ B(i_1,...,i_m),
$$
where $B(i_1,...,i_m)$ is some set of slopes.  We furthermore assume
there is a set $C$ of real numbers such that for each of these
$> N^{m-O(\delta)}$ sets $B(i_1,...,i_m)$, there exists a real number 
$\tau(i_1,...,i_m)$, such that 
\begin{equation} \label{BB}
|B(i_1,...,i_m)\ \triangle\ (C+\tau(i_1,...,i_m))|\ <\ 
|B(i_1,...,i_m)| n^{-O(\delta)}.
\end{equation}
Here, $S \triangle T$ denotes the symmetric difference between $S$
and $T$.
\bigskip

\noindent {\bf Proof of the claim.}  Basically, we just need to show
how these slopes $\theta_i$ link up with the lines in 
(\ref{form1}) and (\ref{form2}); further, we need to explain the
presence of the $\delta$ here, rather than the $\delta'$ appearing
earlier.

Let us first address the issue of the $\delta$ versus of the $\delta'$:  
Since we get to choose $\delta' > 0$ as small as desired relative to
$\eps > 0$, we can just as well rewrite it is $\delta > 0$.

As to the relationship between the $\theta_i$'s above and the $\lambda_j$'s
in (\ref{form1}), we will take 
$$
\{ \theta_1,\ ...,\ \theta_N\}\ =\ \{ \lambda_i \} \cup \{ 1/\lambda_i\}.
$$ 
Then, for $m = 2^{t-2}$ we have that the lines of (\ref{form1}) have
slope of the form $\theta_{i_1} \cdots \theta_{i_m}$.  Furthermore,
the fact that $t > k$ is what will allow us to take $m$ as large as needed.
\hfill $\blacksquare$
\bigskip

Now we combine together pairs of these rich lines -- as discussed in 
subsection \ref{combine_subsection} -- having the same slope,
to produce many other rich lines having slope $1$:  Fix one of the slopes 
$\theta_{i_1}\cdots \theta_{i_m}$ leading to rich lines with
the set of slopes $B(i_1,...,i_m)$.  Applying Lemma \ref{overlap_lemma},
we find that there are at least 
$$
|B(i_1,...,i_m)|^2 n^{-O(\delta)}
$$
ordered pairs
$$
(b,b')\ \in\ B(i_1,...,i_m)\times B(i_1,...,i_m),
$$
such that the line 
$$
y\ =\ x\ +\ (b-b')/\theta_{i_1}\cdots \theta_{i_m}
$$
is $n^{1-O(\delta)}$-rich in the grid $A \times A$.      

From (\ref{BB}), and a little bit of effort, we can easily deduce 
that at least $|B(i_1,...,i_m)|^2 n^{-O(\delta)}$ 
of these pairs $(b,b')$ have the property that there exists
$(c,c') \in C \times C$ satisfying  
$$
(b,b')\ =\ (c + \tau(i_1,...,i_m),c' + \tau(i_1,...,i_m)).
$$
For such pairs, we will have that 
$$
b-b'\ =\ c-c'.
$$

By the pigeonhole principle, there exists at least one
pair (in fact, lots of pairs) $(c,c') \in C \times C$, $c \neq c'$, such that 
at least $N^{m-O(\delta)}$ of the sequences $i_1,...,i_m$ have the 
property that the line 
$$
y\ =\ x\ +\ (c-c')/\theta_{i_1}\cdots \theta_{i_m}
$$
is $n^{1-O(\delta)}$-rich in the grid $A \times A$.  Let us denote
this constant $c - c'$ as $\xi$, so that our rich lines all look like
$$
y\ =\ x\ +\ \xi \phi_{i_1} \cdots \phi_{i_m},\ {\rm where\ } 
\phi_i\ :=\ 1/\theta_i.
$$

By combining together pairs of these lines, as discussed in 
subsection \ref{combine_subsection}, we can form new ones of the form
\begin{equation} \label{repeat_me}
y\ =\ x\ +\ \xi (\phi_{i_1} \cdots \phi_{i_m} - \phi_{j_1} \cdots \phi_{j_m})
\end{equation}
that are rich in the grid.  If we then combine together pairs of {\it those}
lines, we get ones of the form 
\begin{eqnarray} \label{new_lines}
y\ =\ x &+& \xi (\phi_{i_1} \cdots \phi_{i_m} - \phi_{j_1}\cdots \phi_{j_m} 
\nonumber \\
&&\ \ \ \ + \phi_{k_1} \cdots \phi_{k_m} - \phi_{\ell_1} \cdots 
\phi_{\ell_m} ).
\end{eqnarray}
Continuing in this manner, we can generate lines of slope $1$ with 
$y$-intercept equal to $\xi$ times alternating sums of $m$-fold products 
of the $\phi_i$'s; and, at the $t$th iteration, these alternating sums 
have $2^t$ terms.  

\subsection{The sequence $\Theta_i$}

Now we take a digression for a few pages, and define and analyze 
a certain sequence of expressions:  Starting with the set 
$$
\Theta\ :=\ \{ \varphi_i\ :\ i=1,2,... \},
$$
consider the sequence of sets (expressions) 
\begin{equation} \label{thetasequence}
\Theta_1 := \Theta.\Theta - \Theta.\Theta,\ 
\Theta_2 := \Theta_1.\Theta_1 - \Theta_1.\Theta_1,
\end{equation}
and so on.  If we formally expand out the expressions, we will get
sums of the following type:  $\Theta_1$ consists of sums of the type
$$
a_1 a_2 - a_3 a_4,\ a_i \in \Theta,
$$
and $\Theta_2$ consists of the sums
\begin{eqnarray} \label{theta2} 
&& a_1 a_2 a_5 a_6 - a_3 a_4 a_5 a_6 - a_1 a_2 a_7 a_8 + a_3 a_4 a_7 a_8 
\nonumber \\
&& - a_9 a_{10} a_{13} a_{14} + a_9 a_{10} a_{15} a_{16} + a_{11} a_{12}
a_{13} a_{14} - a_{11} a_{12} a_{15} a_{16},
\end{eqnarray}
where again each $a_i \in \Theta$.  We will not bother to write down 
$\Theta_3$!  In general, at the $j$th iteration, the terms in the
alternating sum will involve $4^j$ variables $a_i$, and
the number of terms will be $2^{2^j-1}$.  
\bigskip

Later on, in another subsection, we will show that so long as 
$\delta > 0$ is small enough, upon 
expanding $\Theta_{t - 2}$ into the alternating sum of products of 
variables $a_1,...,a_{4^{t-2}}$, as in (\ref{thetasequence}) and 
(\ref{theta2}), at least  
$$
|\Theta|^{4^{t-2}} n^{-O_t(\delta)}
$$
choices for these $a_i \in \Theta$ will produce a 
$$
\theta = \theta(a_1,...,a_{4^{t-2}}) \in \Theta_{t-2}
$$
so that the line 
\begin{equation} \label{the_lines}
y\ =\ x\ +\ \xi \theta
\end{equation}
is $n^{1 - O_t(\delta)}$-rich in the grid $A \times A$.  We will then 
use Lemma \ref{energy_lemma} to show that this is impossible for 
$t$ large enough and $\delta > 0$ small enough.  The fact that $t > k$,
where $k$ is chosen as large as desired ($k$ is as appears in
subsection \ref{usable_subsection}), will allow us to reach our
contradiction, thereby proving Theorem \ref{main_theorem}.

\subsubsection{A certain inductive claim}

The key fact that we will show and use to accomplish our goal is 
the following.
\bigskip

\noindent {\bf Claim 4.}  Suppose that $g(x_1,...,x_u)$ is some polynomial
in the variables $x_1,...,x_u$, which are to be thought of as taking
on values in the set $\Theta$.  Consider the expansion of 
$$
\Theta_j \Theta_j g(x_1,...,x_u)
$$
into the variables $a_1,...,a_{2 \cdot 4^j}, x_1,...,x_u \in \Theta$.
\footnote{The first $\Theta_j$ is expanded into 
$a_1,...,a_{4^j}$, and the second $\Theta_j$ is expanded into
$a_{4^j+1},...,a_{2 \cdot 4^j}$.}  Suppose that there are at least 
$$
|\Theta|^{2\cdot 4^j + u} n^{-O_{j,u}(\delta)}.
$$
choices for these variables, producing a value 
$$
\gamma\ =\ \gamma(a_1,...,x_u)\ =\ \Theta_j \Theta_j g(x_1,...,x_u)
$$
such that the line 
$$
y\ =\ x\ +\ \xi \gamma
$$
is $n^{1-O_{j,u}(\delta)}$-rich in the grid $A \times A$.  Then, 
there are at least 
$$
|\Theta|^{4^{j+1} + u} n^{-O_{j,u}(\delta)}
$$
choices for the variables 
$$
b_1,...,b_{4^{j+1}}, y_1,...,y_u\ \in\ \Theta
$$
such that the line 
$$
y\ =\ x\ +\ \xi \gamma',\ \gamma' = \gamma'(b_1,...,y_u) \in \Theta_{j+1} 
g(y_1,...,y_u)
$$
is $n^{1-O_{j,u}(\delta)}$-rich in $A \times A$.
\bigskip

\noindent {\bf Proof of the claim.}  Under the hypotheses of the above
claim, the pigeonhole principle implies that for at least 
\begin{equation} \label{choices}
|\Theta|^{4^{j+1}+u} n^{-O_{j,u}(\delta)}
\end{equation}
choices of variables 
$$
b_1,...,b_{2 \cdot 4^j}, c_1,..., c_{2 \cdot 4^j}, x_1,...,x_u\ \in\ \Theta,
$$
we will have that if we let 
$$
\gamma_1\ :=\ \gamma_1(b_1,...,b_{2 \cdot 4^j},x_1,...,x_u)\ \in\ 
\Theta_j \Theta_j g(x_1,...,x_u),
$$
and
$$
\gamma_2\ :=\ \gamma_2(c_1,...,c_{2 \cdot 4^j},x_1,...,x_u)\ \in\ 
\Theta_j \Theta_j g(x_1,...,x_u)
$$
(note that the value of $x_1,...,x_u$ here is the same as for $\gamma_1$),
then both the lines 
$$
y\ =\ x\ +\ \xi \gamma_1\ \ {\rm and\ \ } y\ =\ x\ +\ \xi \gamma_2
$$
are $n^{1-O_{j,u}(\delta)}$-rich in $A \times A$.  Furthermore, by dint of
Lemma \ref{overlap_lemma} and the comments following it, 
we will additionally have that for (\ref{choices})
many choices of the $b_i$'s, $c_i$'s, and $x_i$'s, the pair of lines 
may be combined to produce the new line
$$
y\ =\ x\ +\ \xi (\gamma_1 - \gamma_2),
$$
which will also be $n^{1-O_{j,u}(\delta)}$-rich in $A \times A$.  

This
$$
\gamma_1 - \gamma_2\ =\ (\Theta_j \Theta_j - \Theta_j \Theta_j) 
g(x_1,...,x_u)
$$
has the form $\Theta_{j+1} g(x_1,...,x_u)$.  Clearly this proves the
claim.
\hfill $\blacksquare$
\bigskip

A consequence of this claim, and an easy induction argument (to be
described presently), is that if the number of choices for 
$$
x_1,...,x_{2^Z}\ \in\ \Theta
$$
for which 
\begin{equation} \label{goodlines}
y\ =\ x\ +\ \xi x_1 \cdots x_{2^Z}
\end{equation}
is $n^{1-O_Z(\delta)}$-rich in $A \times A$ is at least 
\begin{equation} \label{thatmany}
|\Theta|^{2^Z} n^{-O_Z(\delta)},
\end{equation}
which it is by the properties of the set $\Theta$ described earlier, 
then there are at least 
$$
|\Theta|^{4^Z} n^{-O_Z(\delta)}
$$
choices for $y_1,...,y_{4^Z} \in \Theta$ such that the line
$$
y\ =\ x + \xi \gamma,\ \gamma = \gamma(y_1,...,y_{4^Z}) \in \Theta_Z
$$
is $n^{1-O_Z(\delta)}$-rich in $A \times A$.  

The way that this is proved is as follows:  First, write the product
$$
x_1 \cdots x_{2^Z}\ =\ (x_1 x_2) (x_3 x_4) \cdots (x_{2^Z-1} x_{2^Z}).
$$
Then, applying the claim to the pair $x_1 x_2$, and then $x_3x_4$, and
so on, we deduce that lots of variable choices make lines
$y = x + \xi \alpha$, $\alpha \in \Theta_1 \cdots \Theta_1$ 
($2^{Z-1}$ copies here), rich in $A \times A$.  Then, the claim is
applied again to the products $\Theta_1 \Theta_1$ (grouped in twos), 
leading to lines $y = x + \xi \beta$, $\beta \in \Theta_2 \cdots \Theta_2$
($2^{Z-2}$ copies here).  Continuing, one reaches lines 
$y = x + \xi \gamma$, $\gamma \in \Theta_Z$, as claimed.
\bigskip

Combining this deduction with Claim 3, we deduce:
\bigskip

\noindent {\bf Claim 5.}  There are at least 
$$
N^{4^{t-2} - O_t(\delta)} 
$$
choices of variables $a_1,...,a_{4^{t-2}} \in \Theta$ such that for 
$\theta = \theta(a_1,...,a_{4^{t-2}}) \in \Theta_{t-2}$, the line 
$$
y\ =\ x\ +\ \xi \theta
$$
is $n^{1-O_t(\delta)}$-rich in $A \times A$.

\subsection{A growth lemma } 

Given a probability measure $f$ supported on a finite set $C$, 
we let $f^*$ denote a certain measure on $CC - CC$ given as follows:
\begin{equation} \label{fstar}
f^*(x)\ :=\ \sum_{c_1 c_2 - c_3 c_4 = x} f(c_1)f(c_2)f(c_3)f(c_4).
\end{equation}

\begin{lemma}  \label{energy_lemma} 
Suppose that $C$ is a finite set of real numbers.  Let $f$ be a measure 
on $C$.  Then,
$$
\max_x f^*(x)\ \ll\ ( \max_x f(x) )^{4/3} (\log |C|)^2.
$$
\end{lemma}

\subsubsection{\bf Proof of Lemma \ref{energy_lemma}}
Let 
$$
M\ :=\ \max_x f(x).
$$

We begin by partitioning the set $C$ into the disjoint sets, some of
which may be empty: 
$$
C\ =\ C_1 \cup C_2 \cup \cdots \cup C_k \cup C_0,
$$
where for $i \geq 1$,
$$
C_i\ :=\ \{ c \in C\ :\ f(c) \in (2^{-i}M, 2^{-i+1} M] \}, 
$$
where $C_0$ is the remaining elements of $C$, and where 
$$
k\ =\ \lfloor 5 \log |C|/\log 2 \rfloor + 1.
$$

We define
$$
f^*_{\alpha,\beta,\gamma,\delta} (x)\ :=\ 
\sum_{c_1 \in C_\alpha, c_2 \in C_\beta, c_3 \in C_\gamma, 
c_4 \in C_\delta \atop c_1 c_2 - c_3 c_4 = x} f(c_1)f(c_2)f(c_3)f(c_4).
$$
We have that 
$$
f^*(x)\ =\ \sum_{0 \leq \alpha,\beta,\gamma,\delta \leq k} 
f^*_{\alpha,\beta,\gamma,\delta}(x).
$$
To prove the theorem, then, all we need to do is get bounds on these
individual terms, and then sum them up. 
\bigskip

First, we can easily bound the total contribution of the terms
where any of the $\alpha,\beta,\gamma,$ or $\delta$ is $0$:  
The contribution of all such terms is clearly bounded from above by 
$$
\ll\ \sum_{x \in CC-CC} M 2^{-5 \log |C|/\log 2}
\ \ll\ |C|^{-1}.
$$
\bigskip

Now we handle the other terms.  First, suppose that 
$1 \leq \alpha,\beta,\gamma,\delta \leq k$.  Then, one easily sees 
from the fact $f$ is a probability measure that 
$$
|C_i|\ \ll\ 2^i M^{-1},\ i=\alpha,\beta,\gamma,\delta.
$$
The size of $f^*_{\alpha,\beta,\gamma,\delta}(x)$ is 
\begin{equation} \label{right_quantity}
\ll\ M^4 2^{-\alpha-\beta-\gamma-\delta} 
|\{ a \in C_\alpha, b \in C_\beta, c \in C_\gamma, 
d \in C_\delta\ :\ ab-cd = x\}|.
\end{equation}
To bound this last factor from above, we will apply Elekes's 
\cite{elekes} idea of using the Szemer\'edi-Trotter incidence theorem 
\cite{szemeredi} to prove sum-product inequalities.  We begin with
the Szemer\'edi-Trotter theorem: 

\begin{theorem}  Suppose that one has $N$ points and $L$ lines in 
the plane.  Then, the number of incidences is bounded from above by
$$
O((N L)^{2/3} + N + L).
$$
\end{theorem}

The way we apply this theorem is as follows:  Consider the family of 
lines
$$
ax + cy\ =\ z,\ {\rm where\ } a \in C_\alpha,\ c \in C_\gamma.
$$
Note that there are $|C_\alpha|\cdot |C_\gamma|$ lines in total.  

Each of these lines intersects the grid $C_\beta \times C_\delta$ in some
number of points (or perhaps no points at all).  The total number 
of incidences $(x,y) \in C_\beta \times C_\delta$ is the right-most
factor of (\ref{right_quantity}).  From the Szemer\'edi-Trotter theorem, 
this number is
\begin{eqnarray*}
&\ll&\ (|C_\alpha| \cdot |C_\beta| \cdot |C_\gamma| \cdot |C_\delta|)^{2/3}
+ |C_\beta|\cdot |C_\delta| + |C_\alpha|\cdot |C_\gamma| \\
&\ll&\ 2^{2(\alpha + \beta + \gamma + \delta)/3} M^{-8/3} + 
2^{\beta + \delta} M^{-2} + 2^{\alpha + \gamma} M^{-2}.
\end{eqnarray*}
The total weight $f(a)f(x)f(c)f(y)$ that each such 
representation $ax + cy = z$ gets is 
$$
\ll\ 2^{-\alpha - \beta - \gamma - \delta} M^4.
$$ 
So, 
$$
f^*_{\alpha,\beta,\gamma,\delta}(z)\ \ll\ 
2^{-(\alpha + \beta + \gamma + \delta)/3} M^{4/3} + 
2^{-\alpha - \gamma} M^2 + 2^{-\beta - \delta} M^2.
$$ 
It follows that for all $z \in CC-CC$,
$$
f^*(z)\ \ll\ |C|^{-1} + M^{4/3} (\log |C|)^2\ \ll\ M^{4/3} (\log |C|)^2.  
$$
The second inequality here comes from the fact that $M \geq |C|^{-1}$,
which follows from the fact that $f$ is a probability measure. 

\subsection{Continuation of the proof}

We now define a sequence of functions by first letting 
$$
f_0(h)\ :=\ \left \{ \begin{array}{rl} 1/N,\ &{\rm if\ }h \in \Theta; \\
0,\ &{\rm if\ } h \not \in \Theta. \end{array}\right.
$$
(Note that $f_0$ is a probability measure.)  Then, we inductively define
$$
f_{i+1}(h)\ :=\ f^*_i(h),
$$   
where $f^*$ is as in (\ref{fstar}).  It is easy to see that these 
$f_i$ are all also probability measures.  

The connection between this function $f$ and our sequence of 
$\Theta_i$ is as follows:  For a given real number $h$ we have that 
$f_j(h)$ is $|\Theta|^{-4^j}$ times the number of choices for
$$
x_1,...,x_{4^j}\ \in\ \Theta
$$
such that 
$$
\theta\ =\ \theta(x_1,...,x_{4^j})\ \in\ \Theta_j
$$
satisfies 
$$
\theta\ =\ h.
$$
As will will see, the upper bound on $f_j(h)$ provided by 
Lemma \ref{energy_lemma} will produce for us a lower bound on the 
number of rich lines in our grid.

Now, Lemma \ref{energy_lemma} implies that for some constant $c > 0$,
if 
$$
t\ \geq\ k\ :=\ c \log(1/\eps),
$$ 
then for all $h$,
$$
f^*_{t-2}(h)\ \leq\ 1/n^5
$$
So, for each real number $h$, there are at most 
$$
n^{-5} |\Theta|^{4^{t-2}} 
$$
choices for $x_1,...,x_{4^{t-2}}\in \Theta$ such that 
$\theta = \theta(x_1,...,x_{4^{t-2}})$ equals $h$.  
Combining this with Claim 5, we quickly deduce that there are 
$n^{5 - O_t(\delta)}$ distinct values of $\theta$ among these rich lines
(of Claim 5).  If $\delta > 0$ is small enough relative to 
$\eps$, then we will see that this number exceeds $n^4$.  

We have now reached a contradiction, since there can be at most 
$n^4$ lines that hit an $n \times n$ grid in at least two points each.
Our theorem is now proved.

\section{Acknowledgements}

We would like to thank Boris Bukh and Jozsef Solymosi for ponting 
out a mistake in an earlier version of the main theorem related 
to an application of the Bourgain-Chang Theorem (we had thought it 
was proved for the reals, when in fact it had only been shown for 
the rationals).  

We would like to thank Boris Bukh for pointing out to us the 
``symmetrization'' argument for passing from grids $A \times B$
to $A \times A$.  During an earlier presentation of our proof, we 
had used a more complicated argument to do this.

We would also like to thank Jozsef Solymosi for pointing out some references
to the works of Sudakov, Szemer\'edi and Vu, which although not essential
for the current paper, was essential for \cite{croot}. 

We would like to thank P. M. Wood for his mentioning that some of our
techniques were similar to those appearing in the Bourgain, Katz and Tao
paper.

We would like to thank T. Tao for mentioning that perhaps a certain
theorem of Bourgain might give a simpler proof of our theorem.

We would like to thank H. Helfgott for pointing out section 2.6 of the
book by Tao and Vu.

And finally, we would like to thank A. Granville for reminding us of
Solymosi's comments regarding the Bourgain-Chang Theorem.

\end{document}